\documentclass[11pt]{article}
\usepackage{mathrsfs}
\usepackage{}
\usepackage{amsfonts}
\usepackage{amssymb}
\usepackage{graphicx}
\usepackage{amsmath}\usepackage[T1]{fontenc}\usepackage{lmodern}
\usepackage{esint}

\def\sqr#1#2{{\vcenter{\vbox{\hrule height.#2pt
              \hbox{\vrule width.#2pt height#1pt \kern#1pt \vrule
width.#2pt}
              \hrule height.#2pt}}}}
\def\signed #1{{\unskip\nobreak\hfil\penalty50
              \hskip2em\hbox{}\nobreak\hfil#1
              \parfillskip=0pt \finalhyphendemerits=0 \par}}
\def\endpf{\signed {$\sqr69$}}
\def\i{\infty}
\def\3n{\negthinspace \negthinspace \negthinspace }
\def\2n{\negthinspace \negthinspace }
\def\1n{\negthinspace }

\def\ms{\medskip}

\def\({\Big (}
\def\){\Big )}
\def\[{\Big[}
\def\]{\Big]}

\def\be{\begin{equation}}
\def\bel{\begin{equation}\label}
\def\ee{\end{equation}}
\def\bea{\begin{eqnarray}}
\def\eea{\end{eqnarray}}
\def\bt{\begin{theorem}}
\def\et{\end{theorem}}
\def\bc{\begin{corollary}}
\def\ec{\end{corollary}}
\def\bl{\begin{lemma}}
\def\el{\end{lemma}}
\def\bp{\begin{proposition}}
\def\ep{\end{proposition}}
\def\br{\begin{remark}}
\def\er{\end{remark}}
\def\ba{\begin{array}}
\def\ea{\end{array}}
\def\bd{\begin{definition}}
\def\ed{\end{definition}}

\newtheorem{lemma}{Lemma}[section]
\newtheorem{remark}{Remark}[section]

\newtheorem{theorem}{Theorem}[section]
\newtheorem{corollary}{Corollary}[section]

\newtheorem{definition}{Definition}[section]
\newtheorem{proposition}{Proposition}[section]

\makeatletter
   
   \@addtoreset{equation}{section}
\makeatother
\allowdisplaybreaks
\begin{document}

\title{\bf Stability of parabolic equations in non-cylindrical domains}

\author{Lingyang Liu\thanks{School of Mathematical Sciences, South China Normal University, Guangzhou, 510631, China. E-mail address: liuly@scnu.edu.cn \ms } }

\date{}

\maketitle

\begin{abstract}
This paper addresses the stability of a class of parabolic equations in non-cylindrical domains. We investigate the $L^\infty$-stability of systems for both nondegenerate and degenerate cases. Unlike in cylindrical domains, solutions to such problems may not exhibit exponential decay. An interesting phenomenon observed is that degeneracy has a positive impact on $L^\infty$-norm estimates for solutions to the system.
\end{abstract}


\section{Introduction and main result}\label{sec1}
Let $l$ be a $C^1$ function defined on $\mathbb R$. For $t\geq0$, set $\Omega_t=(0,l(t))$ and $Q_t=\cup_{0<t<\infty}\Omega_t\times \{t\}$.
Consider the following the parabolic problem:
\begin{equation}\label{e1}
\left\{\begin{array}{ll}
y_{t}-(a(x)y_x)_x=0, & (x,t)\in Q_t,\\[2mm]
y(0,t)=0, \ y(l(t),t)=0, &t\in (0,\infty),\\[2mm]
y(x,0)=y_0(x),&x\in(0,L_0),
\end{array}\right.
\end{equation}
where $y$ is the state variable, $a(x)=x^\alpha$ ($0 \leq \alpha < 1$), and $l(t)=L_0(1+kt)^\gamma$ with $L_0>0$, $k>0$ and $\gamma\in \mathbb R$.

For any $y_0\in L^2(0,L_0)$,  system \eqref{e1} is well-posed in
appropriate functional spaces.

The stability of \eqref{e1} is defined in the following $L^\infty$-sense:
\begin{definition}
System \eqref{e1} is said to be \item[(1)] {exponentially stable} if there exist constants $C > 0$, $\widetilde{C} > 0 $, $t_0 > 0 $ and $\epsilon \geq 1 $ such that the solutions of (\ref{e1}) satisfy
\begin{equation*}
\|y(\cdot, t)\|_{L^\infty(\Omega_t)} \leq C \|y_0\|_{L^2(\Omega_0)} e^{-\widetilde{C}t^\epsilon}, \quad \forall y_0 \in L^2(\Omega_0), \ \forall t \geq t_0.
\end{equation*}
\item[(2)] {subexponentially stable}  if there exist constants $C > 0$, $\widetilde{C} > 0 $, $t_0 > 0 $ and $0 < \sigma < 1 $ such that its solutions satisfy
\begin{equation*}
\|y(\cdot, t)\|_{L^\infty(\Omega_t)} \leq C \|y_0\|_{L^2(\Omega_0)} e^{-\widetilde{C} t^{\sigma}}, \quad \forall y_0 \in L^2(\Omega_0), \ \forall t \geq t_0.
\end{equation*}
\item[(3)] {polynomially stable} if there exist constants $ C > 0 $, $\tau > 0 $ and $t_0 > 0 $  such that its solutions satisfy
\begin{equation*}
\|y(\cdot, t)\|_{L^\infty(\Omega_t)} \leq C \|y_0\|_{L^2(\Omega_0)} (\psi(t))^{-\tau}, \quad \forall y_0 \in L^2(\Omega_0), \ \forall t \geq t_0,
\end{equation*}
where $ \psi(t) $ is a polynomial with respect to $ t $.
\end{definition}
\begin{remark}
As can be seen from the above definitions, the types of stability form a hierarchy of decreasing strength: exponential stability implies subexponential stability, and subexponential stability implies polynomial stability.
\end{remark}

Stability issues for parabolic equations in cylindrical domains have been extensively studied (see, e.g., \cite{BT,CMP,HWW,KP,LR,Nic,Wan} and the rich literature cited therein). However, to our knowledge, very few results exist on the stability of parabolic equations with moving boundaries. In \cite{GLL}, the authors investigated the stability of the degenerate heat equation with a linear moving boundary $l(t)=1+kt$, $k>0$, while \cite{LG} studied the heat equation subject to boundary functions of the form $l(t)=(1+kt)^\gamma$,  $\gamma\in\mathbb R$. Both works focus on the $L^2$-stability of systems.

Our objective is to analyze the $L^\infty$-norm decay properties of solutions to system \eqref{e1}, where the principal coefficient $a(x) = x^\alpha$ degenerates at $x = 0$ for $\alpha>0$. Notably, in the strongly degenerate case $\alpha \geq 1$, solutions to \eqref{e1} may fail to admit a trace on the degenerate boundary. Hence, we are concerned with the $L^\infty$-stability of \eqref{e1} for $0 \leq\alpha<1$. To derive the $L^\infty$-norm estimates, we employ a weighted multiplier
method, which is inspired by the works of \cite{CLG, Liu, SLL}, where the authors established energy estimates and observability inequalities for wave equations in noncylindrical domains. The key to the proof lies in the choice of appropriate weight functions. Furthermore, spectral analysis is used to seek exact solutions in specific cases, and the comparison principle is instrumental in analyzing the decay behavior of solutions.

The main results of this paper are as follows:

\begin{theorem}\label{cor2}
Let $\alpha =0$. The stability of system \eqref{e1}  is characterized as follows:
\begin{itemize}
    \item[(i)] \eqref{e1} is exponentially stable in $L^\infty$-sense for $\gamma\leq0$;
    \item[(ii)] \eqref{e1} is subexponentially stable in $L^\infty$-sense for $0<\gamma<\frac{1}{2}$;
    \item[(iii)]\eqref{e1} is only polynomially stable in $L^\infty$-sense for $\gamma\geq\frac{1}{2}$.
\end{itemize}
\end{theorem}

\begin{remark}
Here, $\gamma=\frac{1}{2}$ serves as a critical value separating subexponential and polynomial stability.
\end{remark}

\begin{theorem}\label{t1.2}
Let $0 < \alpha < 1$. System \eqref{e1} is exponentially stable for $\gamma\leq0$, subexponentially stable for $0<\gamma <\frac{1}{2-\alpha}$, only polynomially stable for $\gamma =\frac{1}{2-\alpha}$, and stable but not exponentially or subexponentially stable for $\gamma >\frac{1}{2-\alpha}$.
\end{theorem}

\begin{remark}
When $\gamma=\frac{1}{2}$, in the nondegenerate case $\alpha = 0$, system \eqref{e1} is only polynomially stable according to the conclusion (iii) of Theorem \ref{cor2}. However, it is subexponentially stable for all $0 < \alpha < 1$, which indicates that degeneracy plays a positive role in the stability of the system.
\end{remark}

\begin{remark}
Given $0 \leq\alpha < 1$, the stability of system \eqref{e1} with a general boundary function $x = S(t)$ can be assessed by the comparison principle.
\end{remark}

The rest of this paper is organized as follows. In Section \ref{sec2}, we analyze the stability of nondegenerate parabolic equations and prove Theorem \ref{cor2}. Section \ref{sec3} is devoted to studying the degenerate parabolic problem and proving Theorem \ref{t1.2}.

\section{Stability of nondegenerate parabolic systems}\label{sec2}

\subsection{$L^\infty$-norm estimate for nondegenerate parabolic equations}

In this section, we investigate the stability of system \eqref{e1} in the nondegenerate case $\alpha = 0$. To this
aim, we first establish a $L^\infty$-norm estimate for a class of nondegenerate parabolic equations. Let $T > 0$. For any function $l\in C^{1}([0,T])$ with $l(0)=L_0$, set  $Q^l_T=\{(x,t):x\in(0,l(t)),t\in(0,T)\}$. Consider the following parabolic problem:
\begin{equation}\label{e2.1}
\left\{\begin{array}{ll}
y_{t}-(ay_x)_x=0, & (x,t)\in Q^l_T,\\[2mm]
y(0,t)=0, \ y(l(t),t)=0, &t\in (0,T),\\[2mm]
y(x,0)=y_0(x),&x\in(0,L_0),
\end{array}\right.
\end{equation}
where $a$ is a continuously differentiable function satisfying
\begin{equation}\label{c1}
a(x, t) \geq \rho,\quad  (x,t)\in Q^l_T,
\end{equation}
for some constant $\rho>0$.

By \cite{LMZ}, for any
$y_0\in L^2(0,L_0)$,
 system \eqref{e2.1} admits a unique solution
\begin{equation*}
y\in C([0, T]; L^2(0,l(t)))\cap L^2(0, T; H_0^1(0,l(t))).
\end{equation*}

\begin{remark}
Solutions to (\ref{e2.1}) possess higher regularity at time $T>0$ due to the regularizing effect.
\end{remark}

As a preliminary, we introduce the extremum principle for (\ref{e2.1}):
\begin{lemma}\label{l2.2}
Assume that $y$ is the classical solution to system (\ref{e2.1}).
Denote the parabolic boundary by $\Gamma_T := \partial{Q^l_T}\setminus[0, l(T)]$. Then
\begin{equation*}
\max\limits_{(x,t)\in \overline{Q^l_T}}
y(x, t) = \max\limits_{(x,t)\in \Gamma_T}y(x, t),
\quad
\min\limits_{(x,t)\in \overline{Q^l_T}}
y(x, t) = \min\limits_{(x,t)\in \Gamma_T}y(x, t).
\end{equation*}
\end{lemma}

Now, we prove the following $L^\infty$-norm estimate for system \eqref{e2.1}.
\begin{proposition}\label{l2.3}
The solution $y$ of \eqref{e2.1} satisfies the following estimates:
\begin{equation}\label{2.2}
\|y(\cdot, T)\|_{L^\infty(0,l(T))}\leq{\|y_0\|_{L^2(0,L_0)}}\bigg(\int_{0}^{T}\frac{\rho e^{\rho\alpha(t)}}{l(t)}dt\bigg)^{-\frac{1}{2}},
\end{equation}
where $\alpha(t)= \int_{0}^{t} \frac{1}{l^2(\tau)}d\tau$.
\end{proposition}
\emph{ Proof of Proposition \ref{l2.3}.}
The proof is divided into two parts.

{\bf Step 1.} First, we define the maximum function $M(t):=\max\limits_{x\in[0,l(t)]}|y(x, t)|$. By Lemma \ref{l2.2},
 $M(t)$ is nonincreasing with respect to $t$, i.e.,
\begin{equation}\label{e3.1}
M( t_1)\geq M( t_2) ,\quad 0<t_1\leq t_2.
\end{equation}

On the other hand, by the Newton-Leibniz formula and the Cauchy-Schwarz inequality, it is easy to see that
\begin{equation*}
 |y(x, t)|=\left|\int_{0}^{x}y_v(v,t)dv\right|\leq\int_{0}^{l(t)}|y_v(v,t)|dv\leq \sqrt{l(t)}\Big(\int_{0}^{l(t)}|y_v(v,t)|^2dv\Big)^{\frac{1}{2}},\quad x \in[0, l(t)].
\end{equation*}
Squaring both sides, we find that
\begin{equation}\label{M2(t)}
M^2(t)\leq l(t)\int_{0}^{l(t)}|y_x(x,t)|^2dx.
\end{equation}

{\bf  Step 2.} Next, multiply both sides of the first equation in (\ref{e2.1}) by $e^{s\alpha(t)}y$ (where $s>0$ and the weight function $\alpha(t)$ will be specified later), integrate over $(0, l(t))$, perform integration by parts, and use the boundary conditions $y(0,t)=y(l(t), t) = 0$ to get
\begin{equation*}
\frac{d}{2dt}\int_{0}^{l(t)}e^{s\alpha(t)}y^2dx-\frac{1}{2}\int_{0}^{l(t)}s\alpha_te^{s\alpha(t)}y^2dx=-\int_{0}^{l(t)}e^{s\alpha(t)}ay^2_xdx.
\end{equation*}
By the assumption \eqref{c1}, it follows that
\begin{equation}\label{e3.2}
\frac{d}{2dt}\int_{0}^{l(t)}e^{s\alpha(t)}y^2dx\leq-\rho\int_{0}^{l(t)}e^{s\alpha(t)}y^2_xdx+\frac{1}{2}\int_{0}^{l(t)}s\alpha_te^{s\alpha(t)}y^2dx.
\end{equation}
Noticing that
\begin{equation*}\label{3.0}
\begin{array}{rl}
\displaystyle\int_{0}^{l(t)}y^2(x,t)dx
&\displaystyle\!= \int_{0}^{l(t)}\Big(\int_{0}^{x}y_v(v,t)dv\Big)^2dx \\[3mm]
&\displaystyle\leq \int_{0}^{l(t)}x\int_{0}^{x}y^2_v(v,t)dvdx \\[3mm]
&\displaystyle\leq l^2(t)\int_{0}^{l(t)}y^2_x(x,t)dx,
\end{array}
\end{equation*}
we substitute the choices $\alpha_t = l^{-2}(t) $, $s=\rho$ into \eqref{e3.2} to derive
\begin{equation}\label{e3.5}
\frac{d}{dt}\int_{0}^{l(t)}e^{\rho\alpha(t)}y^2dx\leq-\rho\int_{0}^{l(t)}e^{\rho\alpha(t)}y^2_xdx.
\end{equation}
This, together with  \eqref{M2(t)}, indicates
\begin{equation*}
\frac{d}{dt}\int_{0}^{l(t)}e^{\rho\alpha(t)}y^2dx\leq\frac{-\rho e^{\rho\alpha(t)} M^2(t)}{l(t)}.
\end{equation*}
Integrating the above inequality over $[0,T]$ and applying \eqref{e3.1}, we obtain
\begin{equation*}
\int_{0}^{l(T)}e^{\rho\alpha(T)}y^2(x,T)dx-\int_{0}^{L_0}e^{\rho\alpha(0)}y^2_0(x)dx\leq-\rho\int_{0}^{T}\frac{e^{\rho\alpha(t)}M^2(t)}{l(t)}dt\leq-\rho M^2(T)\int_{0}^{T}\frac{e^{\rho\alpha(t)}}{l(t)}dt.
\end{equation*}
With the initial condition $\alpha(0) = 0$, the above inequality leads to
\begin{equation*}
-\int_{0}^{L_0}y^2_0(x)dx\leq-\rho M^2(T)\int_{0}^{T}\frac{e^{\rho\alpha(t)}}{l(t)}dt.
\end{equation*}
Consequently,
\begin{equation*}\label{3.3}
M(T)\leq{\|y_0\|_{L^2(0,L_0)}}{}\bigg(\int_{0}^{T}\frac{\rho e^{\rho\alpha(t)}}{l(t)}dt\bigg)^{-\frac{1}{2}}.
\end{equation*}
This completes the proof of Proposition \ref{l2.3}. \endpf

\begin{remark}
An $L^2$-norm estimate for system \eqref{e2.1} can be obtained from \eqref{e3.5}, which implies that
\begin{equation*}
\frac{d}{dt}\int_{0}^{l(t)}e^{\rho\alpha(t)}y^2dx\leq0.
\end{equation*}
Hence,
\begin{equation*}
\int_0^{l(T)}e^{\rho\alpha(T)}y^2(x,T)dx\leq\int_0^{L_0}e^{\rho\alpha(0)}y^2_0(x)dx.
\end{equation*}
That is,
\begin{equation*}
\displaystyle\|y(T)\|_{L^2(0,l(T))}\leq e^{\frac{{\rho}}{2}(\alpha(0)-\alpha(T))}\|y_0\|_{L^2(0,L_0)}=e^{-\frac{{\rho}}{2}\int^T_0\frac{1}{l^2(t)}dt}\|y_0\|_{L^2(0,L_0)},
\end{equation*}
which is applicable in analyzing the $L^2$-stability of nondegenerate parabolic equations.
\end{remark}

\subsection{Proof of Theorem \ref{cor2}}

Now, we are in a position to give a proof of Theorem \ref{cor2}.

Let $l(t)=L_0(1+kt)^\gamma$, with $k > 0$ and $\gamma\in \mathbb R$. We derive the $L^\infty$-norm decay rate of solutions to the heat equation:
\begin{equation}\label{2.1}
\left\{\begin{array}{ll}
y_{t}-ay_{xx}=0, & (x,t)\in Q_t,\\[2mm]
y(0,t)=0, \ y(l(t),t)=0, &t\in (0,\infty),\\[2mm]
y(x,0)=y_0(x),&x\in(0,L_0),
\end{array}\right.
\end{equation}
where the diffusion coefficient $a$ is a positive constant.

Recalling that $\alpha(t) = \int_{0}^{t} l^{-2}(\tau) d\tau$, we compute that
\begin{equation}\label{3.4}
\alpha(t)=\left\{\begin{array}{ll}
\left[L^2_0(1-2\gamma)k\right]^{-1}\left[(1+kt)^{1-2\gamma}-1\right],&\gamma \neq\frac{1}{2},\\[3mm]
\displaystyle\frac{1}{L^2_0k}\ln(1+kt),&\gamma =\frac{1}{2}.
\end{array}\right.
\end{equation}

$(i)$ When $\gamma\leq0$,  \eqref{3.4} yields
\begin{equation*}
\begin{array}{ll}
\alpha(t) \geq \left[L_0^2(1-2\gamma)\right]^{-1}t, \quad  t \geq 0.
\end{array}
\end{equation*}
Moreover, notice that $l(t)\leq L_0$ for all $t \geq 0$. It follows that
\begin{equation*}
\begin{array}{rl}
\displaystyle\frac{e^{\rho\alpha(t)}}{l(t)}\geq L^{-1}_0e^{\rho\alpha(t)}&\displaystyle\!\!\!\geq L^{-1}_0e^{\rho\left[L^2_0(1-2\gamma)\right]^{-1}t},\quad  t \geq 0,
\end{array}
\end{equation*}
which leads to
\begin{equation}\label{s2}
\begin{array}{rl}
\displaystyle\int_{0}^{T}\frac{e^{\rho\alpha(t)}}{l(t)}dt\geq L^{-1}_0\int_{0}^{T} e^{\rho\left[L^2_0(1-2\gamma)\right]^{-1}t}dt=\frac{L_0\left(1-2\gamma\right)}{\rho}\left[ e^{\rho\left[L^2_0(1-2\gamma)\right]^{-1}T}-1\right].
\end{array}
\end{equation}
Substituting \eqref{s2} into \eqref{2.2}, we deduce that there exists $t_0 \geq 0$ such that
\begin{equation*}
\begin{array}{rl}
\displaystyle\|y(\cdot, T)\|_{L^\infty(0,l(T))}\leq\sqrt{2}\left[L_0\left(1-2\gamma\right)\right]^{-\frac{1}{2}}e^{-\frac{\rho}{2}\left[L^2_0(1-2\gamma)\right]^{-1}T}\|y_0\|_{L^2(0,L_0)}, \quad T\geq t_0,
\end{array}
\end{equation*}
which means that system \eqref{2.1} is exponentially stable in the $L^\infty$ sense for $\gamma\leq0$.

$(ii)$ When $0<\gamma < \frac{1}{2}$, it holds that $l(t)\geq L_0$ for all $t \geq 0 $. Thus,
\begin{equation*}
\begin{array}{rl}
\displaystyle\int_{0}^{T}\frac{e^{\rho\alpha(t)}}{l(t)}dt\geq L_0\int_{0}^{T}\frac{e^{\rho\alpha(t)}}{l^2(t)}dt&\displaystyle\!\!\!=L_0\int_{0}^{T}e^{\rho\alpha(t)}d\alpha(t)=\frac{L_0}{\rho}\left[e^{\rho\alpha(T)}-1\right].
\end{array}
\end{equation*}
Since $0<1-2\gamma<1$, it follows from \eqref{3.4} that for sufficiently large $T>0$,
\begin{equation*}
\alpha(T)\geq \frac{1}{2}\left[L^2_0(1-2\gamma)k^{2\gamma}\right]^{-1}T^{1-2\gamma}.
\end{equation*}
Therefore,
\begin{equation*}
\begin{array}{rl}
\displaystyle\int_{0}^{T}\frac{e^{\rho\alpha(t)}}{l(t)}dt\geq \frac{L_0}{2\rho}e^{\rho\alpha(T)}\geq \frac{L_0}{2\rho}e^{\frac{\rho}{2}\left[L^2_0(1-2\gamma)k^{2\gamma}\right]^{-1}T^{1-2\gamma}}.
\end{array}
\end{equation*}
Consequently, \eqref{2.2} indicates that
\begin{equation*}
\begin{array}{rl}
\displaystyle\|y(\cdot, T)\|_{L^\infty(0,l(T))}\leq\sqrt{2}\left(L_0\right)^{-\frac{1}{2}}e^{-\frac{\rho}{4}\left[L^2_0(1-2\gamma)k^{2\gamma}\right]^{-1}T^{1-2\gamma}}\|y_0\|_{L^2(0,L_0)},
\end{array}
\end{equation*}
which implies that system \eqref{2.1} is subexponentially stable in the $L^\infty$ sense for $0<\gamma<\frac{1}{2}$.

$(iii)$ {\bf Case 1.} When $\gamma = \frac{1}{2}$, we use the second expression for $\alpha(t)$ in \eqref{3.4} to get
\begin{equation*}
\begin{array}{rl}
\displaystyle\int_{0}^{T}\frac{e^{\rho\alpha(t)}}{l(t)}dt= \frac{1}{L_0}\int_{0}^{T}(1+kt)^{\frac{\rho}{L^2_0k}-\frac{1}{2}}dt\!\!&\displaystyle=\frac{2L_0}{2\rho+L^2_0k}\left[(1+kT)^{\frac{\rho}{L^2_0k}+\frac{1}{2}}-1\right].
\end{array}
\end{equation*}
This implies that for sufficiently large $T>0$,
\begin{equation*}
\begin{array}{rl}
\displaystyle\int_{0}^{T}\frac{e^{\rho\alpha(t)}}{l(t)}dt\geq\displaystyle \frac{L_0}{2\rho+L^2_0k}(1+kT)^{\frac{\rho}{L^2_0k}+\frac{1}{2}}.
\end{array}
\end{equation*}
Substituting this into \eqref{2.2}, we obtain
\begin{equation}\label{4.3}
\begin{array}{rl}
\displaystyle\|y(\cdot, T)\|_{L^\infty(0,l(T))}\leq\left(\frac{L_0\rho}{2\rho+L^2_0k}\right)^{-\frac{1}{2}}(1+kT)^{-\frac{\rho}{2L^2_0k}-\frac{1}{4}}\|y_0\|_{L^2(0,L_0)}.
\end{array}
\end{equation}
From \eqref{4.3}, we conclude that system \eqref{2.1} is polynomially stable in the $L^\infty$ sense for $\gamma = \frac{1}{2}$.

On the other hand, we demonstrate that system \eqref{2.1} is only polynomially stable for $\gamma = \frac{1}{2}$.

Let $y$ be the solution to system \eqref{2.1}. We introduce new variables
$$\eta=\frac{x}{L_0(1 + kt)^{\frac{1}{2}}},\quad \theta = \ln(1+kt).$$
Put $v(\eta, \theta)=y(x,t)$. Then $v$ is the solution to
\begin{equation}\label{3.10}
\left\{\begin{array}{ll}
v_{\theta}-{{k}^{-1}L^{-2}_0}av_{\eta\eta}-\frac{1}{2}\eta v_\eta=0, & (\eta,\theta)\in Q_\Theta,\\[2mm]
v(0,\theta)=0, \ v(1,\theta)=0, &\theta\in (0,\Theta),\\[2mm]
v(\eta,0)=y_0(L_0\eta),&x\in(0,1),
\end{array}\right.
\end{equation}
where $Q_\Theta=(0,1)\times(0,\Theta)$ with $\Theta=\ln(1+kT)$.

Assume the solution $v$ has the form $v(\eta, \theta)=X(\eta)T(\theta)$. Substituting it into the equation \eqref{3.10} and applying the method of separation of variables, we obtain \begin{equation*}
T_{\theta}+\mu T=0,
\end{equation*}
and the associated eigenvalue problem
\begin{equation}\label{4.5}
\left\{\begin{array}{ll}
{{k}^{-1}L^{-2}_0}aX_{\eta\eta}+\frac{1}{2}\eta X_\eta+\mu X=0, \\[2mm]
X(0)=0, \ X(1)=0.
\end{array}\right.
\end{equation}
Let $p(\eta)=e^{\frac{kL_0^2\eta^2}{4a}}$ and $AX=-\frac{{k}^{-1}L^{-2}_0a}{p}(pX_\eta)_\eta$. We define the weighted spaces
\begin{equation*}
\begin{array}{ll}
L^{2}_{p}(0,1)=\{X:(0,1)\rightarrow \mathbb R|\int^1_0p(\eta)X^2d\eta<\infty\}, \\[2mm]
H^{1}_{p}(0,1)=\{X\in L^{2}_{p}(0,1) |X_\eta\in L^{2}_{p}(0,1)\}, \\[2mm]
H^{0,1}_{p}(0,1)=\{X\in H^{1}_{p}(0,1) |X(0)=X(1)=0\},\\[2mm]
D(A)=\{X\in H^{0,1}_{p}(0,1) |AX\in L^{2}_{p}(0,1)\}.
\end{array}
\end{equation*}
$A: D(A)\subset L^{2}_{p}(0,1)\rightarrow L^{2}_{p}(0,1)$ is an unbounded operator. By an argument similar to that in \cite{MZ}, we conclude that the eigenvalues $\mu_n$ of \eqref{4.5} are positive, satisfy $0 < \mu_1 \leq \mu_2 \leq \cdots \leq \mu_n \leq \cdots$, and tend to infinity as $n \to \infty$. Additionally, the corresponding eigenvectors $\{ \phi_n \}_{n\in \mathbb N^\ast}$ form a complete orthogonal system in $L^{2}_{p}(0,1)$.

Taking $v(\eta,0)= \phi_n(\eta)$, the corresponding solution to \eqref{3.10} will be
\begin{equation*}
v(\eta, \theta)=e^{-\mu_n\theta}\phi_n(\eta).
\end{equation*}
Returning to the original variables, the solution $y$ takes the form
\begin{equation*}
y(x, t)=(1+kt)^{-\mu_n}\phi_n\bigg(\frac{x}{L_0(1 + kt)^{\frac{1}{2}}}\bigg).
\end{equation*}
Consequently, its $L^\infty$-norm satisfies
\begin{equation*}
\|y(\cdot, t)\|_{L^\infty(0,l(t))}=(1+kt)^{-\mu_n}\|\phi_n\|_{L^\infty(0,1)},
\end{equation*}
which implies that system (\ref{2.1}) is only polynomially stable when $\gamma= \frac{1}{2}$.

{\bf Case 2.} When $\gamma > \frac{1}{2}$, we prove that system (\ref{2.1}) is also polynomially stable. Consider the following Cauchy problem
\begin{equation}\label{4.1}
\left\{\begin{array}{ll}
z_{t}-az_{xx}=0,&  (x,t)\in {\mathbb R}\times(0,T),\\[2mm]
z(x,0)=z_0(x),&x\in\mathbb R.
\end{array}\right.
\end{equation}
It is well known that problem \eqref{4.1} admits a unique solution of the form
\begin{equation*}
z(x,t) =(4\pi a t)^{-\frac{1}{2}}\int_{\mathbb R}z_0(v)e^{-\frac{|x-v|^2}{4at}}dv.
\end{equation*}
Without loss of generality, we suppose $y_0\in H^1_0(0,L_0)$. Let $z_0\geq0$ (or $z_0\leq0$) be continuously differentiable and absolutely integrable on $\mathbb{R}$, and satisfy
\begin{equation*}
z_0(x)=
\left\{\begin{array}{ll}
\|y_0\|_\infty\ (\text{or} -\|y_0\|_\infty),&x\in(0,L_0),\\[1mm]
0,&\mathbb R\backslash(0,2L_0).
\end{array}\right. \qquad
\end{equation*}
By the comparison principle, it follows that
\begin{equation}\label{4.2}
\sup\limits_{x\in (0,l(t))}|y(\cdot,t)|\leq \sup\limits_{x\in\mathbb R}|z(\cdot,t)|\leq (4\pi at)^{-\frac{1}{2}}\|z_0\|_{L^1(0,2L_0)},
\end{equation}
where $y$ is the solution to system \eqref{2.1} with the initial data $y_0$.

\eqref{4.2} means that solutions to system \eqref{2.1} decay at a rate no slower than $t^{-\frac{1}{2}}$ in the $L^\infty$ sense.

Finally, we demonstrate that system (\ref{2.1}) lacks subexponential stability for $\gamma>\frac{1}{2}$. Put
$$
L(t) = L_0(1 + kt)^{\frac{1}{2}}, \quad \widehat{\Omega}_t = (0, L(t)).
$$
Let ${y}_1$ be the solution to \eqref{2.1} associated with initial data $\phi_1$ and the boundary $L$, satisfying
\begin{equation}\label{4.4}
\| y_1(\cdot, t) \|_{L^\infty(\widehat{\Omega}_t)} =(1+kt)^{-\mu_1}\|\phi_1\|_{L^\infty(0,1)}.
\end{equation}
Take $y_0 \geq \max\{\phi_1^+,\phi_1^- \}$, where $\phi_1^+$ and $\phi_1^-$ denote the positive and negative parts of $\phi_1$, respectively.
Let $y$ be the solution to \eqref{2.1} corresponding to $y_0$ and $l(t)=(1+kt)^{\gamma}$, $\gamma>\frac{1}{2}$. Then ${y}|_{\widehat{\Omega}_t \times (0, \infty)}$ is the
solution to the system
\begin{equation*}
\begin{cases}
 \hat{y}_t - a\hat{y}_{xx} = 0, & \text{in } \widehat{\Omega}_t \times (0, \infty), \\
 \hat{y}(0, t) =0,\hat{y}(L(t), t) = \mu(t), & \text{in } (0, \infty), \\
  \hat{y}(x, 0) = y_0(x), & \text{in } \widehat{\Omega}_0,
\end{cases}
\end{equation*}
where $\mu(t) ={y}(L(t), t)$.
By the comparison principle, we obtain
\begin{equation}\label{4.8}
\| {y}(\cdot, t) \|_{L^\infty(\Omega_t)} \geq \| \hat{y}(\cdot, t) \|_{L^\infty(\widehat{\Omega}_t)} \geq \| y_1(\cdot, t) \|_{L^\infty(\widehat{\Omega}_t))}.
\end{equation}
\eqref{4.2}, together with \eqref{4.4} and \eqref{4.8}, indicates that system \eqref{2.1} is only polynomially stable for $\gamma >\frac{1}{2}$.

The proof of Theorem \ref{cor2} is complete. \endpf

\section{Stability of degenerate parabolic systems}\label{sec3}

\subsection{Degenerate parabolic equations in cylindrical domains}
In order to prove Theorem \ref{t1.2}, we first study a class of degenerate parabolic equations defined in cylindrical domains.

Given $0<\alpha<1$, define
$
a(x)=x^\alpha,\ \forall x\geq0.
$
For $T > 0$, set $Q_T :=(0, 1)\times(0, T )$. Consider the initial-boundary value problem
\begin{equation}\label{e4.1}
\left\{\begin{array}{ll}
u_{t}-(a(x,t)u_{x})_x+b(x,t)u_x=0,&  (x,t)\in Q_T,\\[2mm]
u(0,t)=u(1,t)=0,&  t\in (0,T),\\[2mm]
u(x,0)=u_0(x),&x\in(0,1),
\end{array}\right.
\end{equation}
where ${a}(x,t)=p(t)a(x)$, $b(x,t)=q(t)x$, $p\in  C^1(\overline{Q_T })$ is positive in $\overline{Q_T}$, $\frac{1}{p} \frac{\partial p}{\partial t} \in L^\infty(Q_T)$ and $q\in L^{\infty}(Q_T)$.

The following definition is given in \cite{Wang}.
\begin{definition}
Define $\mathscr B$ to be the closure of the set $C_0^\infty(Q_T)$ with respect to
the norm
\begin{equation*}
\|u\|_{\mathscr B} =\left(
\iint_{Q_T}
a(x, t)\left(|u|^2 + | u_x|^2\right)dxdt
\right)^{1/2}
,\ u \in \mathscr B.
\end{equation*}
\end{definition}
Then, for any $u_0\in L^2(0,1)$, there exists a unique weak solution $u \in L^\infty((0, T ); L^2(0,1)) \cap\mathscr B$ to the problem (\ref{e4.1}).
\begin{remark}\label{t}
Let $\Sigma = \{0,1\}\times(0, T)$. If $u \in \mathscr B$, then $u|_\Sigma= 0$ in the trace sense, while there is no trace on $\left(\{0\}\times(0, T)\right)$ in general when $\alpha\geq1$.
\end{remark}
By applying the parabolic regularization method, one can prove the following extremum principle:
\begin{lemma}\label{l4.0}
If $u$ satisfies equation (\ref{e4.1}), then
\begin{equation*}
\sup\limits_{(x,t)\in {Q_T}}
u(x, t) \leq \sup\limits_{(x,t)\in \partial_p Q_T}u(x, t),
\quad
\inf\limits_{(x,t)\in {Q_T}}
u(x, t) \geq \inf\limits_{(x,t)\in \partial_p Q_T}u(x, t),
\end{equation*}
where $\partial_p Q_T:= \partial Q_T\setminus \left( [0,1] \times \{T\} \right)$.
\end{lemma}
As an application of Lemma \ref{l4.0}, the following comparison principle holds.
\begin{lemma}\label{l4.1}
Suppose $v$ and $w$ satisfy equation (\ref{e4.1}) with $v|_{\partial_p Q_T} \leq w|_{\partial_p Q_T}$. Then
\begin{equation*}
v(x, t) \leq w(x, t),\quad \forall (x,t)\in Q_T.
\end{equation*}
\end{lemma}

We recall the the following Hardy-type inequality (see \cite{CMV}).
\begin{lemma}\label{l4.2}
Let $0 \leq \alpha < 1$. Then, for all locally absolutely continuous functions $u$ on
$(0, 1)$ satisfying
\begin{equation*}
\underset{x \to 0^+}{z(x)\rightarrow 0}  \ \text{and}\ \int_0^{1}x^\alpha u_x^2dx < \infty,
\end{equation*}
the following inequality holds:
\begin{equation}\label{hd}
\int_0^{1}x^{\alpha-2}u^2dx\leq \frac{4}{(1-\alpha)^2}\int_0^{1}x^\alpha u_x^2dx.
\end{equation}
\end{lemma}

Now, we establish the $L^\infty$-norm estimate for solutions to system \eqref{e4.1}.
\begin{proposition}\label{p2.3}
If a function $\beta: [t_0,T] \to \mathbb{R}$ fulfills
\begin{equation}\label{4.11}
4\left(|\beta_t(t)| + |q(t)|\right) \leq(1-\alpha)^2 p(t),\quad t \in [t_0,T],
\end{equation}
then solutions of \eqref{e4.1} satisfy the estimate
\begin{equation}\label{4.6}
\sup\limits_{x\in[0,1]}|u(x, T)|\leq\|u(t_0)\|_{L^2(0,1)}\left((1-\alpha)e^{-\beta(t_0)}\int_{t_0}^{T}e^{\beta(t)}p(t) dt\right)^{-\frac{1}{2}}.
\end{equation}
\end{proposition}

\emph{Proof.}  The proof is divided into two parts.

{\bf Step 1.} For any positive integer number $n$, choose $a_n$, $b_n$, $u^{(n)}_0\in C^\infty(\overline{Q_T})$ satisfying
\begin{equation}\label{e4.2}
\begin{array}{ll}
 a(x,t)+\frac{1}{n}\leq a_{n}(x,t)\leq a(x,t)+\frac{2}{n},\quad \|b_{n,x}\|_{L^{\infty}(0,1)}\leq |q(t)|,\quad \|u^{(n)}_0\|_{L^2(0,1)}\leq \|u_0\|_{L^2(0,1)}.
\end{array}
\end{equation}

Consider the problem
\begin{equation}\label{e4.3}
\left\{\begin{array}{ll}
u^{(n)}_{t}-(a_{n}(x,t)u^{(n)}_{x})_x+b_{n}(x,t)u^{(n)}_x=0,&  (x,t)\in Q_T,\\[2mm]
u^{(n)}(0,t)=u^{(n)}(1,t)=0,&  t\in (0,T),\\[2mm]
u^{(n)}(x,0)=u_0(x),&x\in(0,1).
\end{array}\right.
\end{equation}
According to the classical theory of parabolic equations, the problem (\ref{e4.3}) admits a unique classical solution $u^{(n)}$. First, we define $S^{(n)}(t):=\sup\limits_{x\in[0,1]}|u^{(n)}(x, t)|$. By the maximum principle,
 $S^{(n)}(t)$ is nonincreasing with respect to $t$, i.e.,
\begin{equation}\label{3.1}
S^{(n)}( t_1)\geq S^{(n)}( t_2) ,\quad 0<t_1\leq t_2.
\end{equation}

On the other hand, by the Newton-Leibniz formula and the Cauchy-Schwarz inequality, it is easy to see that for any $x \in[0, 1]$,
\begin{equation*}
 |u^{(n)}(x, t)|=\left|\int_{x}^{1}u^{(n)}_v(v,t)dv\right|\leq\int_{0}^{1}|u^{(n)}_v(v,t)|dv\leq \bigg(\int_{0}^{1}v^{-\alpha} dv\bigg)^{\frac{1}{2}}\bigg(\int_{0}^{1}v^\alpha |u^{(n)}_v|^2dv\bigg)^{\frac{1}{2}}.
\end{equation*}
Thus, for any $0<\alpha<1$,
\begin{equation}\label{M(t)}
|S^{(n)}(t)|^2\leq \frac{1}{1-\alpha}\int_{0}^{1}x^\alpha|u^{(n)}_x|^2dx.
\end{equation}

{\bf  Step 2.}
Multiply both sides of the first equation in (\ref{e4.3}) by $e^{\beta(t)}u^{(n)}$, integrate over $(0,1)$, then apply integration by parts to get
\begin{equation*}
\frac{d}{2dt}\int_{0}^{1}e^{\beta(t)}|u^{(n)}|^2dx-\frac{1}{2}\int_{0}^{1}\beta_te^{\beta(t)}|u^{(n)}|^2dx=-\int_{0}^{1}e^{\beta(t)}a_{n}|u^{(n)}_x|^2dx+\frac{1}{2}\int_{0}^{1}e^{\beta(t)} b_{n,x}| u^{(n)}|^2dx.
\end{equation*}
By \eqref{e4.2}, it follows that
\begin{equation}\label{3.2}
\frac{d}{2dt}\int_{0}^{1}e^{\beta(t)}|u^{(n)}|^2dx\leq-\int_{0}^{1}e^{\beta(t)}a|u^{(n)}_x|^2dx+\frac{1}{2}\int_{0}^{1}e^{\beta(t)}\left(|\beta_t|+|q(t)|\right)| u^{(n)}|^2dx.
\end{equation}
Since $\int_{0}^{1}| u^{(n)}|^2dx\leq C_\alpha\int_{0}^{1}x^{\alpha}|u^{(n)}_x|^2dx$ (see \eqref{hd}), by choosing $C_\alpha\left(|\beta_t| + |q(t)|\right) \leq p(t)$ in \eqref{3.2}, we derive
\begin{equation*}\label{3.5}
\frac{d}{dt}\int_{0}^{1}e^{\beta(t)}|u^{(n)}|^2dx\leq-\int_{0}^{1}e^{\beta(t)}a|u^{(n)}_x|^2dx.
\end{equation*}
This, together with  \eqref{M(t)}, indicates
\begin{equation*}
\frac{d}{dt}\int_{0}^{1}e^{\beta(t)}|u^{(n)}|^2dx\leq-(1-\alpha)e^{\beta(t)}p(t)|S^{(n)}(t)|^2.
\end{equation*}
Integrating the above inequality over $[t_0,T]$ and applying \eqref{3.1}, we obtain
\begin{equation*}
\begin{array}{ll}
&\displaystyle\int_{0}^{1}e^{\beta(T)}|u^{(n)}(T)|^2dx-\int_{0}^{1}e^{\beta(t_0)}|u^{(n)}(t_0)|^2dx\\[2mm]
&\displaystyle\leq-(1-\alpha)\int_{t_0}^{T}e^{\beta(t)}p(t)|S^{(n)}(t)|^2dt\\[2mm]
&\displaystyle\leq-(1-\alpha)|S^{(n)}(T)|^2\int_{t_0}^{T}e^{\beta(t)}p(t) dt.
\end{array}
\end{equation*}
The above inequality leads to
\begin{equation*}
\int_{0}^{1}e^{\beta(t_0)}|u^{(n)}(t_0)|^2dx\geq(1-\alpha)|S^{(n)}(T)|^2\int_{t_0}^{T}e^{\beta(t)}p(t) dt,
\end{equation*}
which yields
\begin{equation*}\label{3.3}
S(T):=\sup\limits_{x\in[0,1]}|u(x, T)|\leq\|u(t_0)\|_{L^2(0,1)}\left((1-\alpha)e^{-\beta(t_0)}\int_{t_0}^{T}e^{\beta(t)}p(t) dt\right)^{-\frac{1}{2}}.
\end{equation*}
This completes the proof of Proposition \ref{p2.3}. \endpf

\begin{remark}
If $q<0$, \eqref{4.6} remains valid provided that
\begin{equation*}
4|\beta_t(t)|  \leq(1-\alpha)^2 p(t),\quad t \in [t_0,T].
\end{equation*}
\end{remark}

\subsection{Proof of Theorem \ref{t1.2}}
Now, we are in a position to give a proof of Theorem \ref{t1.2}.

Let $l(t)=L_0(1 + kt)^{\gamma}$, with $k > 0$ and $\gamma \in \mathbb{R}$. We study the stability of the following degenerate
linear system:
\begin{equation}\label{eq}
\left\{\begin{array}{ll}
y_{t}-(a(\xi)y_{\xi})_\xi=0,&  (\xi,t)\in Q_t,\\[2mm]
y(0,t)=y(l(t),t)=0,&  t\in(0,\infty),\\[2mm]
y(\xi,0)=y_0(\xi),&\xi\in(0,L_0),
\end{array}\right.
\end{equation}
where $a(\xi)=\xi^\alpha$$(0<\alpha<1)$ and $Q_t=\{(\xi,t):\xi\in(0,l(t)),t\in(0,\infty)\}$.

Put $\xi=l(t)x$ and $u(x,t)=y(\xi,t)$. Then \eqref{eq} can be transformed into a system of the form
\begin{equation}\label{4.0}
\left\{\begin{array}{ll}
u_{t}-p(t)(x^\alpha u_{x})_x+q(t)xu_x=0,&  (x,t)\in Q_T,\\[2mm]
u(0,t)=u(1,t)=0,&  t\in (0,T),\\[2mm]
u(x,0)=u_0(x),&x\in(0,1),
\end{array}\right.
\end{equation}
where
\begin{equation}\label{p}
p(t)=l^{\alpha-2}(t)\quad \text{and} \quad q(t)=-k\gamma(1+kt)^{-1}.
\end{equation}
Note that
\begin{equation}\label{n}
\|y(\cdot,t)\|_{L^\infty(0,l(t))}=\|u(\cdot,t)\|_{L^\infty(0,1)},\quad \forall t\geq0.
\end{equation}
Hence, proving Theorem \ref{t1.2} reduces to establishing the analogous results for system \eqref{4.0}.

$(1)$ When $\gamma <\frac{1}{2-\alpha}$, there exist $t_0>0$ and $\varepsilon>0$ such that
\begin{equation*}
\varepsilon(1+kt)^{\gamma(\alpha-2)}+4|q(t)|\leq(1-\alpha)^2 p(t),\quad \forall t\geq t_0,
\end{equation*}
where $p$ and $q$ are given in \eqref{p}.

Hence, the choice $\beta(t)=\frac{\varepsilon}{4k[\gamma(\alpha-2)+1]}(1+kt)^{\gamma(\alpha-2)+1}$ satisfies \eqref{4.11}. Moreover,
\begin{equation*}
\int_{t_0}^{T}e^{\beta(t)}p(t) dt=4L_0^{\alpha-2}{\varepsilon}^{-1}\int_{t_0}^{T}e^{\beta(t)} d{\beta(t)}.
\end{equation*}
Substituting this into \eqref{4.6}, we obtain
\begin{equation*}
\sup\limits_{x\in[0,1]}|u(x, T)|\leq C(\|u(t_0)\|_{L^2(0,1)},\alpha,\beta(t_0),\varepsilon,L_0)\|u_0\|_{L^2(0,1)}e^{-\frac{1}{2}\beta(T)}.
\end{equation*}
Since ${\gamma(\alpha-2)+1}\geq1$ for $\gamma\leq0$ and $0<{\gamma(\alpha-2)+1}<1$ for $0<\gamma <\frac{1}{2-\alpha}$,
it follows that \eqref{4.0} is exponentially stable when $\gamma\leq0$ and subexponentially stable when $0<\gamma <\frac{1}{2-\alpha}$.

In particular, when $\gamma =\frac{1}{2}$, system \eqref{4.0} exhibits subexponential stability.

$(2)$ When $\gamma =\frac{1}{2-\alpha}$, \eqref{4.0} can be written as
\begin{equation}\label{4.10}
\left\{\begin{array}{ll}
u_{t}=\frac{1}{1+kt}\left[L_0^{\alpha-2}(x^\alpha u_{x})_x+\frac{k}{2-\alpha}xu_x\right],&  (x,t)\in Q_T,\\[2mm]
u(0,t)=u(1,t)=0,&  t\in (0,T),\\[2mm]
u(x,0)=u_0(x),&x\in(0,1).
\end{array}\right.
\end{equation}
Let $\beta(t)=\frac{\varepsilon}{k}\ln(1+kt)$ with $4\varepsilon\leq(1-\alpha)^2L_0^{\alpha-2}$. By an argument analogous to Step 2 in the proof of Proposition \ref{p2.3}, we get
\begin{equation*}
\sup\limits_{x\in[0,1]}|u(x, T)|\leq C(\alpha,\varepsilon,k)\|u_0\|_{L^2(0,1)}(1+kT)^{-\frac{\varepsilon}{2k}}.
\end{equation*}
This implies that the decay of solutions to system \eqref{4.0} is at least polynomial.

Next, we consider the boundary value
problem
\begin{equation}\label{4.7}
\left\{\begin{array}{ll}
-\left[L_0^{\alpha-2}(x^\alpha u'(x))'+\frac{k}{2-\alpha}xu'(x)\right]=\lambda u(x),&  x\in (0,1), 0<\alpha<1, \lambda\in\mathbb R,\\[2mm]
u(0)=u(1)=0.
\end{array}\right.
\end{equation}
The first equation in (\ref{4.7}) may also be written as
\begin{equation}\label{4.9}
(P(x)u'(x))'+\lambda K(x) u(x)=0,
\end{equation}
where
\begin{equation*}
P(x)=L_0^{\alpha-2}K(x)x^\alpha,\quad K(x)=e^{L_0^{2-\alpha}\frac{k}{(2-\alpha)^2}x^{2-\alpha}}.
\end{equation*}
For any $\alpha \in[0, 1)$, we define the space of square-integrable functions on the interval $[0,1]$ weighted by $K(x)$ as $L_{2,K}$, equipped with the inner product
$$
(u_1,u_2)_K =\int^1_0K(x)u_1(x)u_2(x)dx,\quad u_1,u_2\in L_{2,K}.
$$
Additionally, define the inner product in $C^1_0[0,1]$ as follows:
\begin{equation}\label{H}
(u_1,u_2)_H =\int^1_0P(x)u'_1(x)u'_2(x)dx,\quad u_1,u_2\in C^1_0[0,1].
\end{equation}
The norm induced by (\ref{H}) is denoted by $\|\cdot\|_H$, and the completion of $C^1_0[0,1]$ with respect to this norm, denoted by $H^{0,1}_{P}$, is a Hilbert space. Let
\begin{equation*}
J(u)=\frac{(u,u)_H}{(u,u)_K},\quad \forall u\in H^{0,1}_{P}.
\end{equation*}
The first eigenvalue and the corresponding eigenfunction of problem \eqref{4.7} are given by
\begin{equation*}
\begin{array}{ll}
\lambda_1 = \inf\limits_{\substack{u \in H^{0,1}_{P} \\ u \neq 0}} J(u),\quad J(\phi_1)=\lambda_1.
\end{array}
\end{equation*}
Furthermore, there exist a sequence $(\lambda_n)_{n\in \mathbb N^\ast}$ of real numbers with $\lambda_n > 0$ and
$\lambda_n \rightarrow \infty$, and the corresponding
eigenfunctions $(\phi_n)_{n\in \mathbb N^\ast}$ in $H^{0,1}_{P}$ such that \eqref{4.9} holds.

For a fixed positive integer $n$, let $u_0(x)=\phi_n(x)$. Assume that the solution to \eqref{4.10} takes the form $u(x,t) = \phi_n(x)\psi_n(t)$. Multiplying both sides of the first equation in \eqref{4.10} by $u$ and integrating over $(0,1)$, we derive
\begin{equation*}
\frac{d}{2dt}(u,u)_{L^2(0, 1)}=-\frac{\lambda_n}{1+kt}(u,u)_{L^2(0, 1)}.
\end{equation*}
It follows that
\begin{equation*}
\|u(t)\|^2_{L^2(0, 1)}=(1+kt)^{-\frac{2\lambda_n}{k}}\|\phi_n\|^2_{L^2(0, 1)}.
\end{equation*}
Furthermore,
\begin{equation*}
\|u(t)\|_{L^\infty(0, 1)}\geq\|u(t)\|_{L^2(0, 1)}=(1+kt)^{-\frac{\lambda_n}{k}}\|\phi_n\|_{L^2(0, 1)}.
\end{equation*}
This implies that \eqref{4.0} is only polynomially stable for $\gamma =\frac{1}{2-\alpha}$.

In light of \eqref{n}, returning to the original variables, the above conclusions continue to hold for system \eqref{eq}.

$(3)$ When $\gamma >\frac{1}{2-\alpha}$, we prove that \eqref{eq} is at most polynomially stable by the comparison
principle. It holds that
\begin{equation*}
l(t)=L_0(1+kt)^{\gamma}> L(t):=L_0(1+kt)^{\frac{1}{2-\alpha}}.
\end{equation*}
Take $y_0 \geq\max\{\phi_n^+,\phi_n^-\}$, where $\phi_n^+$ and $\phi_n^-$ are the positive and negative parts of the eigenfunction $\phi_n$ for problem \eqref{4.7}.

Denote by $y$ the solution to \eqref{eq} associated with the initial data $y_0$ and the moving boundary $l$. It is clear from Lemma \ref{l4.0} that $y\geq0$ in $Q_t$. Consider the parabolic system:
\begin{equation*}
\left\{\begin{array}{ll}
\hat{y}_{t}-(a(x)\hat{y}_{x})_x=0,&  (x,t)\in Q_t,\\[2mm]
\hat{y}(0,t)=0, \hat{y}(L(t),t)=\mu(t),&  t\in(0,\infty),\\[2mm]
\hat{y}(x,0)=y_0(x),&x\in(0,L_0),
\end{array}\right.
\end{equation*}
where $\mu(t)=y(L(t),t)\geq0$. By the comparison principle (see Lemma \ref{l4.1}), we have
\begin{equation*}
\|y(\cdot, t)\|_{L^\infty(0,l(t))} \geq\|\hat{y}(\cdot, t)\|_{L^\infty(0,L(t))} \geq \|{y}_n(\cdot, t)\|_{L^\infty(0,L(t))},\quad \forall t\geq0,
\end{equation*}
where ${y}_n$ denotes the solution to \eqref{eq} subject to the initial data $\phi_n$ and the moving boundary $L$.

Since $\|{y}_n(\cdot, t)\|_{L^\infty(0,L(t))}\geq(1+kt)^{-\frac{\lambda_n}{k}}\|\phi_n\|_{L^2(0,1)}$, solutions to system \eqref{eq} decay at a rate no faster than polynomial.
The proof of Theorem \ref{t1.2} is complete.
\endpf

\end{document}